\documentclass[12pt]{article}
\usepackage{amssymb}

\def\aa{{\cal A}}
\def\bb{{\mathfrak b}}
\def\dd{{\mathfrak d}}
\def\cc{{\mathfrak c}}
\def\pp{{\mathfrak p}}

\def\dom{{\rm dom}}

\def\swap{{\rm swap}}
\def\proof{\par\noindent Proof\par\noindent}
\def\poset{{\mathbb P}}
\def\qposet{{\widetilde{\mathbb Q}}}
\def\rposet{{\mathbb Q}}
\def\eposet{{\mathbb E}}
\def\forces{{| \kern -2pt \vdash}}
\def\force{\forces}
\def\res{\upharpoonright}
\def\qed{\par\noindent QED\par}

\def\rmand{\mbox{ and }}

\newtheorem{theorem}{Theorem}
\newtheorem{lemma}[theorem]{Lemma}

\newtheorem{conjecture}[theorem]{Conjecture}

\begin{document}

\begin{center}
{\large A MAD Q-set}
\end{center}
\begin{flushright}
Arnold W. Miller\footnote{
Thanks to the Fields Institute for Research in Mathematical Sciences
at the University of Toronto for their support during the time this paper
was written and to Juris Steprans who directed the special program in set
theory and analysis.
\par Mathematics Subject Classification 2000: 03E35
} 
\end{flushright}

\begin{quote}
 \par\centerline{Abstract}\par
A MAD (maximal almost disjoint) family is an infinite
subset $\aa$ of the infinite subsets of
$\omega=\{0,1,2,\ldots,\}$ such that any two elements of
$\aa$ intersect in a finite set and every infinite subset
of $\omega$ meets some element of $\aa$ in an infinite set.
A Q-set is an uncountable set of reals such that every 
subset is a relative $G_\delta$ set.
It is shown that it is relatively consistent with
ZFC that there exists a MAD family which is also a Q-set in
the topology in inherits a subset of $P(\omega)=2^{\omega}$. 
\end{quote}

In this paper we answer a question of Hrusak by showing that it is consistent
that there exists a maximal almost disjoint family  
$\aa\subseteq [\omega]^\omega$ which also a Q-set.  A topological space is a
Q-set iff every subset is a $G_\delta$.  His reason for asking this question
was because in a certain argument involving a topological space 
$\Psi({\cal A})$ built from a  MAD family it would have been helpful to assume
that a MAD family cannot  be a Q-set. Szeptycki \cite{szept} contains some
results on vanDouwen's $\Psi$ and also on Q-sets.

Our construction is similar to that in Fleissner and Miller \cite{FM} where a
Q-set is obtained which is concentrated on the rationals.  In Judah and Shelah
\cite{JS} it is shown consistent to have a Q-set while at the same time
$\bb=\dd=\omega_1$. Their Q-set forcing has the Sack's property.  Their forcing
is also used in Nowik and Weiss \cite{nowik} to construct a Q-set with certain
properties and also Gruenhage and Koszmider \cite{gruen} to construct a
topological space with certain properties. In our model as in \cite{FM} we have
that $\dd=\cc=\omega_2$ and $\bb=\omega_1$.  

In Dow \cite{dow} and Brendle \cite{brendle} a type of Q-set forcing is used
which preserves towers (so $\pp=\omega_1$) 
which generalizes Hechler dominating real forcing, and
$\bb=\dd=\cc$.

\begin{theorem}
It is relative consistent with ZFC, that there exists a MAD family
$\aa\subseteq [\omega]^\omega$ which also a Q-set.
\end{theorem}
\proof

We begin by forcing a generic MAD family and then we iterate
our Q-set forcing to make the generic MAD into a Q-set. The difficulty
is to ensure the family stays maximal.  

\bigskip
Define: Let $\poset$ be the usual poset for forcing a MAD family:
\par\noindent $(p,q)\in \poset$ iff 
\begin{enumerate}
 \item $p:F\to 2^{N}$ for some finite $F\subseteq\omega_1$ and $N<\omega$
 (write $F=\dom(p)$ and $N=N_p$)
 \item $q$ a partial function from a subset of $[F]^2$ into $N$
 \item if $q(\alpha,\beta)=n$, then for every $i$ with $n\leq i<N$ either
 $p(\alpha)(i)=0$ or $p(\beta)(i)=0$.
\end{enumerate}

The uniformity $N$ of lengths in condition(2) is not strictly necessary but
it will be convenient and would occur on a dense set anyway.

\bigskip
Define: 
\par\noindent $(p_1,q_1)\leq (p_2,q_2)$ iff 
\begin{enumerate}
 \item $\dom(p_1)\supseteq \dom(p_2)$
 \item $p_1(\alpha)\supseteq p_2(\alpha)$ for all $\alpha\in dom(p_2)$
 \item $q_1\supseteq q_2$ 
\end{enumerate}

\noindent Intuitively, we are describing a family 
$\{a_\alpha\subseteq\omega:\alpha<\omega_1\}$ as follows:

\begin{enumerate}
 \item $p(\alpha)=s$ means  ($i\in a_\alpha$ iff $s(i)=1$) for $i<|s|$
 \item $q(\alpha,\beta)=n$ promises that $a_\alpha\cap a_\beta\subseteq n$ 
\end{enumerate}

Note that $(p_1,q_1)$ and $(p_2,q_2)$ are compatible iff there exists
$p_3\leq p_1,p_2$ such that $(p_3,q_1\cup q_2)$ is in $\poset$.

This forcing is due to Hechler \cite{hechler}.  For $G$ 
$\poset$-generic over $M$ define
     $$x_\alpha^G=\bigcup\{p(\alpha): \exists q \;(p,q)\in G\}$$
And let $X=\{x_\alpha^G:\alpha<\omega_1\}$ and
let $\aa=\{a_\alpha\subseteq \omega:\alpha<\omega\}$ where
each $x_\alpha$ is the characteristic function of $a_\alpha$, i.e.
$$a_\alpha=\{n:x_\alpha(n)=1\}$$

The following lemma is due to Hechler.
\begin{lemma} $\poset$ is ccc.  If $G$ is $\poset$-generic over
$M$, then in $M[G]$ the set $\aa$ is a maximal almost disjoint family of
infinite subsets of $\omega$.
\end{lemma}

We will in a sense need to reprove this lemma since we will show
that after our new version of Q-set forcing our generic
family still remains
a maximal almost disjoint family.
The idea of the argument is that given a name $\tau$
for some infinite subset of $\omega$, we find an $\alpha$ 
which is not involved with deciding $n\in\tau$ for any $n$.
Then we get a contradiction by swapping the value of
$x_\alpha(n)=0$ to $x_\alpha'(n)=1$ while still forcing $n\in\tau$.
In the usual Q-set forcing while the condition forcing $n\in\tau$ doesn't
directly talk about $x_\alpha$, it may decide that $[s]\subseteq U_n$
where the other condition says $x_\alpha\notin U_n$.  These conditions
may become inconsistent when we change to $x_\alpha'$ because it
might be that $s\subseteq x_\alpha'$ even though $s$ is not
a subset of $x_\alpha$.

\begin{center}
A new Q-set forcing
\end{center}

The following is to motivate our definition of $\poset*\qposet$.  It
would be the definition of the new Q-set forcing in the model
$M[G]$ where $G$ is $\poset$-generic.

\bigskip
Define. For $x\in 2^\omega$,$s\in 2^{<\omega},k<\omega$
  $$\swap(x,s,k)=\{y\in 2^\omega: s\subseteq y, 
               |\{i\geq |s|:y(i)\not=x(i)\}|\leq k\}$$
Note that $\swap(x,s,k)$ is a countable closed subset of
$[s]$. It contains $x$ if $s\subseteq x$. Also 
$\swap(x,\langle\rangle,0)=\{x\}$.

Suppose we are given $X\subseteq 2^\omega$ such that for all $x\not=y\in X$
there are infinitely many $n$ with $x(n)\not=y(n)$.  For $Y\subseteq X$
define $\rposet(X,Y)$ as follows:

$r\in\rposet(X,Y)$ iff $r$ is a finite subset of
$$\{(n,s):n<\omega, s\in 2^{<\omega}\}\cup\{(n,(x,t,k)): x\in Y:
t\in 2^{<\omega},\; n,k<\omega\}$$
subject to the condition:

\bigskip
\centerline{ if $(n,s)\in r$ and $(n,(x,t,k))\in r$, then
$[s]\cap \swap(x,t,k)=\emptyset$.}
\bigskip

The ordering is by inclusion $r_1\leq r_2$ iff $r_1\supseteq r_2$.
The meaning of these conditions is
\begin{enumerate}
 \item $(n,s)$ means ``$[s]\subseteq U_n$''
 \item $(n,(x,t,k))$ means ``$\swap(x,t,k)\cap U_n=\emptyset$''
\end{enumerate}
Now suppose $G$ is $\rposet(X,Y)$-generic over a model $N$.  Define
$$U_n^G=\bigcup\{[s]: \exists r\in G \;(n,s)\in r\}$$
An easy genericity argument shows that
$$X\cap\bigcap_{n<\omega}U_n^G= X\setminus Y$$
To see this suppose $y\in Y$ and $r$ any condition, let $n$
be sufficiently large so as to not appear in $r$ at all.  
Then let $r'=r\cup\{(n,(y,\langle\rangle,0)\}$ and note that
$$r'\forces y\notin U_n$$  
On the otherhand
let $y\in X\setminus Y$, $r$ any condition, and $n<\omega$ arbitrary.
Since $y$ is infinitely often different from any  element of $X$ mentioned
in $r$ (they must come from $Y$),  we can find $l<\omega$ so that 
$$[y\res l]\cap \swap(x,s,k)=\emptyset$$
for any $(n,(x,s,k))\in r$.  Now we let $r'=r\cup\{(n,y\res l)\}$ then
$$r'\forces y\in U_n$$

Next we describe the ordering $\poset *\qposet$ which is a basic
building block of our iteration.  
If $G$ is $\poset$-generic over $M$ then 
$\qposet^G$ is essentially the same as 
$\rposet(X,X)$.

\bigskip
Define.
\par\noindent $((p,q),r)\in \poset*\qposet$
iff 
\begin{enumerate}
 \item $(p,q)\in\poset$

 \item $r$ is a finite subset of the union of
 $$\{(n,t):n<N_p,\;t\in 2^{<N_p}\}$$ and 
 $$\{(n,(\alpha,s,k)): \alpha\in \dom(p), s\in 2^{<N_p}, \; n, k <N_p\}$$

 \item if  $(n,(\alpha,s,k))\in r$ and $(n,t)\in r$, then either 
  $s$ and $t$ are incomparable or $s\subseteq t$ and 
  $$|\{i: |s|\leq i<|t|, \;\;t(i)\not=p(\alpha)(i)\}|>k$$
\end{enumerate} 
Condition (3) guarantees that for any $x\in 2^\omega$ such
that $x\supseteq p(\alpha)$ that 
$$\swap(x,s,k)\cap [t]=\emptyset$$

\noindent The ordering is given by

$((p_1,q_1),r_1)\leq ((p_2,q_2),r_2)$ iff
$(p_1,q_1)\leq ((p_2,q_2)$ and $r_1\supseteq r_2$.

\noindent Note that 
$((p_1,q_1),r_1)$ and $((p_2,q_2),r_2)$ are compatible iff 
there exists $p_3\leq p_1,p_2$ such that
$((p_3,q_1\cup q_2),r_1\cup r_2)$ is a condition.

\begin{center}
The $\omega_2$ iteration.  
\end{center}

Our iteration can be described as a suborder of the product
$$\poset\times \sum_{\alpha<\omega_2}\eposet$$
Where $\eposet$ is the set of all finite subsets of
$$\{(n,t):n<\omega,t\in 2^{<\omega}\}\cup
\{(n,(\alpha,s,k)): \alpha\in \omega_1,
s\in 2^{<\omega}, n, k <\omega\}$$
and $\sum_{\alpha<\omega_2}\eposet$ is the set of all
$r:\omega_2\to \eposet$ such that $r(\delta)$ is trivial (ie. the
empty set) for all but finitely many $\delta$.

By induction on $\beta\leq\omega_2$ define 
$$\poset_\beta\subseteq \poset\times\sum_{\alpha<\beta}\eposet$$
as follows: 

Define. $\poset_0=\poset$, 

\bigskip
\noindent  Suppose that we have
defined $\poset_\beta$ and we are also given a $\poset_\beta$
name $Y_\beta$
for a subset of $\omega_1$, ie.
$$\forces_\beta Y_\beta\subseteq \omega_1$$

Define. $((p,q),r)\in \poset_{\beta+1}$ iff
\begin{enumerate}
 \item $((p,q),r\res\beta)\in \poset_{\beta}$,
 \item $((p,q),r(\beta))\in \poset*\qposet$ 
 \item $((p,q),r\res\beta)\forces_\beta \alpha\in Y_\beta$ 
  whenever $(n,(\alpha,s,k))\in r(\beta)$ for some $n,s,k,\alpha$
\end{enumerate}

For limit ordinals $\lambda\leq\omega_2$ we define 
$((p,q),r)\in \poset_\lambda$
iff for all $\beta<\lambda$ we have $((p,q),r\res\beta)\in \poset_\beta$
and for all but finitely many $\beta<\lambda$ we have that $r(\beta)$ is
the trivial condition (i.e. empty set).

Since the iteration of ccc forcing is ccc all of these forcings are ccc.  
To see this directly we can argue as follows: 
Standard arguments using $\Delta$ systems show that 
$\poset_\beta$ has precalibre $\omega_1$, ie. 
any $\omega_1$ sequence of conditions contain an $\omega_1$
subsequence which is centered. Start with
$((p_\alpha,q_\alpha),r_\alpha)\in\poset_\beta$ for $\alpha<\omega_1$.
We can find an uncountable $\Sigma\subseteq \omega_1$
and finite sets $F$ and $H$ and $N<\omega$ so that
\begin{enumerate}
\item $N_\alpha=N$ for all $\alpha\in\Sigma$,
\item $\dom(p_\alpha)\cap\dom(p_\beta)=F$ for $\alpha\not=\beta\in\Sigma$,
\item $\dom(r_\alpha)\cap\dom(r_\beta)=H$ for $\alpha\not=\beta\in\Sigma$,
\item $p_\alpha\res F$ are all the
same for $\alpha\in\Sigma$,
\item $q_\alpha\res [F]^2$ are all the
same for $\alpha\in\Sigma$, and
\item $r_\alpha\res H$ are all the
same with respect to $\{(n,s):n<\omega,s\in 2^{<\omega}\}$ 
for $\alpha\in\Sigma$.
\end{enumerate}
Then any two (or even finite subset) of them are compatible.

\bigskip

Assuming that the ground model satisfies the GCH  
by the usual book keeping argument we can 
arrange things
so that for any $Y\subseteq \omega_1$ which appears in
$M[G_{\omega_2}]$ there will be a name for it in the list 
$Y_\alpha$ for some $\alpha<\omega_2$.  The simplest way
to do this is to take 
$$\{\langle Z^\alpha_\beta:\beta<\omega_1\rangle:\alpha<\omega_2\}$$
which lists all $\omega_1$ sequences of
countable subsets of $\poset\times \sum_{\alpha<\omega_2}\eposet$ with
$\omega_2$ repetitions and then define
$$Y_\alpha=\{\langle p, \check{\beta}\rangle: \beta<\omega_1,
p\in Z^\alpha_\beta\cap 
\poset_\alpha\}$$

If we define 
$$x_\alpha=\bigcup\{s\in 2^{<\omega}:\exists ((p,q),r)\in G,\;\;
s=p(\alpha)\}
\rmand X=\{x_\alpha:\alpha<\omega_1\}$$

Then $X$ will be the characteristic functions of an
almost disjoint family
$\aa=\{a_\alpha:\alpha<\omega_1\}$.  Furthermore if
we define the open sets 
$$U_n^\beta=\cup\{[s]: \exists ((p,q),r)\in G\;\; (n,s)\in r(\beta)\}$$
then by the usual genericity argument
$$ \bigcap_{n<\omega}U_n^\beta\cap X=\{x_\alpha:\alpha\notin Y_\beta^G\}$$
and so $X$ will be a $Q$-set.  

\bigskip
The nontrivial part of our argument is to prove that $\aa$ remains 
a maximal almost disjoint family.  So let $\tau$ be a name for a
counterexample, ie. suppose
$$((p_0,q_0),r_0)\forces \tau\in[\omega]^\omega 
\rmand \forall \alpha<\omega_1\; \tau\cap a_\alpha \mbox{ is finite }$$
Let $\Sigma\subseteq\poset_{\omega_2}$ be a countable set of
conditions extending $((p_0,q_0),r_0)$ 
such that for any $n\in\omega\;\;$ $\Sigma$ contains a maximal antichain
beneath $((p_0,q_0),r_0)$ which decides $n\in\tau$. 
Let $\alpha_0<\omega_1$ be any ordinal not mentioned in any condition from
$\Sigma$. We show $a_{\alpha_0}\cap \tau$ is infinite.

Suppose for contradiction that we have $((p_1,q_1),r_1)\leq ((p_0,q_0),r_0)$, 
and $N_1<\omega$ such that
$$((p_1,q_1),r_1)\forces \tau\cap a_{\alpha_0}\subseteq N_1$$
Without loss of generality we may assume that $N_1=N_{p_1}$.
By tacking on strings of zeros to the conditions in $p_1$ we may assume
that every integer occurring in $r$ is bounded by $N_1-2$
(and not just as required by $N_1$). Let 
       $$F=\{\beta:\{\alpha_0,\beta\}\in\dom(q_1)\}$$ 
Define $r'\supseteq r_1$ as follows:

for each $\delta$ 
$$r'(\delta)=r_1(\delta)\cup \{(n,(\alpha_0,t',k+1): 
  (n,(\alpha_0,t,k))\in r_1(\delta),t'\in A_{\delta,n}, t'\supseteq t\}$$
where
$$A_{\delta,n}=\{t'\in 2^{N_1-1}: t'\mbox{ is incomparable with all $s$ such 
that }
 (n,s)\in r_1(\delta)\}$$

Note that $((p_1,q_1),r')$ is a valid condition because $\alpha_0$ is 
forced into $Y_\delta$ and $t'$ incomparable with all $s$ which might
be a problem.  Let $G$ be a generic filter containing 
$((p_1,q_1),r')$.  Since $\tau^G$ is almost disjoint from each $a_\beta^G$
and infinite, there exists some $n_0\in\tau^G$ with $n_0>N_1$ and
$n_0\notin a_\beta^G$ for all $\beta\in F$. Let
$((p_2,q_2),r_2)\in \Sigma\cap G$ be so that
$$((p_2,q_2),r_2)\forces n_0\in \tau$$
Since it is from $\Sigma$ it does not mention $\alpha_0$.

Let $((p^*,q^*),r^*)\in G$ be stronger than both $((p_1,q_1),r')$ and
$(p_2,q_2),r_2)$ and such that $N^*>n_0$.  Note that
$((p^*,q_1\cup q_2),r'\cup r_2)$ is a valid condition. Any $\gamma$ that
needs to be forced into some $Y_\beta$ is already forced in by either
$((p_1,q_1),r'\res\beta)$ or $(p_2,q_2),r_2\res\beta)$. 

If $p^*(\alpha_0)(n_0)=1$ then we already have a contradiction and there is
nothing to prove.  So assume not, and define $p'$ to be exactly the same as
$p^*$ except $p'(\alpha_0)(n_0)=1$.

\bigskip
\noindent {\bf Claim}.
$((p',q_1\cup q_2),r\cup r_2)$ is a valid condition, extending
both $((p_1,q_1),r_1)$ and $(p_2,q_2),r_2)$.

Proof: Note that we have dropped the extra conditions from $r'$, these
were put there just to prove this Claim.  The
fact that $p'$ extends both $p_1$ and $p_2$ uses that $n_0>N_{p_1}=N_1$ and
$\alpha_0$ is not in the domain of $p_2$. Similarly since
$q_2$ does not mention $\alpha_0$, so if 
$\{\alpha_0,\beta\}\in\dom(q_1\cup q_2)$, then $\beta\in F$ and we
know that $p^*(\beta)(n)=0$ for each $\beta\in F$.  So making
$p'(\alpha_0)(n)=1$ does not violate any promises of disjointness
made in $q_1\cup q_2$.
So we have that $(p',q_1\cup q_2)\in\poset$.  

Now fix
$\delta$ and we must check that
$$((p',q_1\cup q_2),r_1(\delta)\cup r_2(\delta))\in \poset*\qposet$$
We need to check condition (3)

(3) if  $(n,(\alpha,s,k)),(n,t)\in r_1(\delta)\cup r_2(\delta)$ 
 then either $s$ and $t$ are 
 incomparable or $s\subseteq t$ and 
 $$|\{i: |s|\leq i<|t|,\; t(i)\not=p'(\alpha)(i)\}|>k$$

Suppose it fails.
It can only fail if the $\alpha=\alpha_0$ and
since $r_2$ does not mention $\alpha_0$ it must be
that $(n,(\alpha_0,s,k))\in r_1(\delta)$ and $(n,t)\in r_2(\delta)$.
Also it must be that $s$ and $t$ are comparable with $s\subseteq t$ but
$$|\{i: |s|\leq i<|t|, \;t(i)\not=p'(\alpha_0)(i)\}|\leq k$$ 

Note also that $|t|>n_0>N_1$ because otherwise
 $$\{i: |s|\leq i<|t|, \;t(i)\not=p'(\alpha_0)(i)\}=
  \{i: |s|\leq i<|t|, \; t(i)\not=p^*(\alpha_0)(i)\}$$
But then
 $$|\{i: |s|\leq i<|t|, \;\;t(i)\not=p^*(\alpha_0)(i)\}|> k$$
because 
$((p^*,q_1\cup q_2),r_1(\delta)\cup r_2(\delta))\in \poset*\qposet$.

Now let $t'=t\res (N_1-1)$. 

\medskip
\noindent Case 1.  $t'$ is comparable with
some  $s'$ such that $(n,s')\in r_1(\delta)$.

Recall that every integer occurring in $r_1(\delta)$ is bounded by
$N_1-1$. So it must be that $s'\subseteq t'$ but intuitively
this is easy because $r_1(\delta)$ is already asserting
$[s']\subseteq U_n^\delta$ and this implies $[t]\subseteq U_n^\delta$.
More formally, 
$s'\subseteq t'$
and therefore $s'$ and $s$ are both initial strings of $t'$ and so 
comparable, but then we know:
 $$|\{i: |s|\leq i<|s'|<N_1 , s'(i)\not=p_1(\alpha_0)(i)\}|>k$$
But this is still true for $p'$ since we have not changed it 
below $N_1$. 

\medskip
\noindent Case 2. $t'\in A_{\delta,n}$ and 
so we added $(\alpha_0,t',k+1)$ to
$r'(\delta)$.

But remember $((p^*,q_1\cup q_2),r'\cup r_2)$ is a valid 
condition, which means that
$$|\{i: N_1\leq i<|t|, t(i)\not=p^*(\alpha_0))(i)\}| > k+1$$
but $k<N_1-2$ and $p^*(\alpha_0)$ agrees with $p'(\alpha_0)$ except
at exactly one coordinate so
$$|\{i: |s|<N\leq i<|t|, t(i)\not=p'(\alpha_0))(i)\}|>k$$

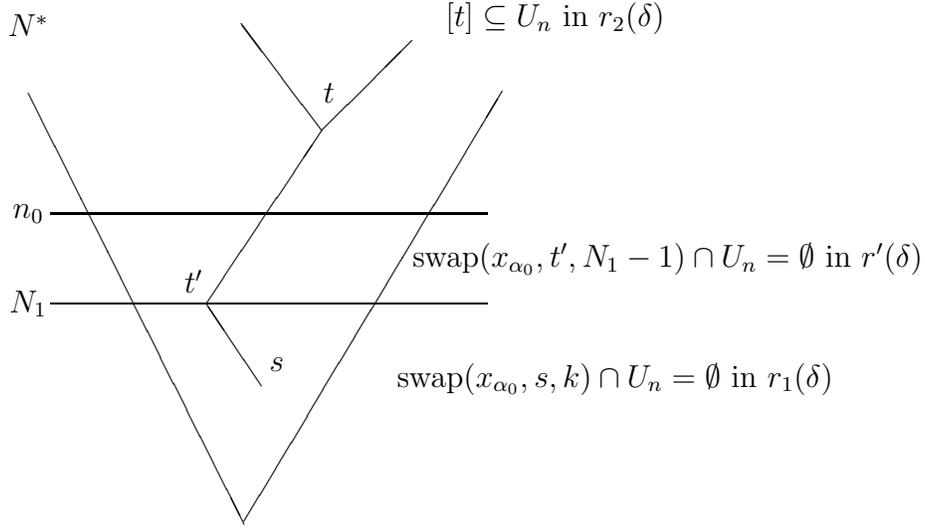
\begin{figure}
\unitlength=1.00mm
\begin{picture}(89.00,79.00)
\put(10.00,70.00){\line(1,-2){28.67}}
\put(38.67,13.00){\line(3,5){34.33}}
\put(13.00,42.00){\line(1,0){58.00}}
\put(10.00,42.00){\makebox(0,0)[cc]{$N_1$}}
\put(41.00,31.00){\line(-2,3){7.33}}
\put(33.67,42.00){\line(2,3){15.33}}
\put(49.00,65.00){\line(-3,4){10.67}}
\put(49.00,65.00){\line(1,1){12.00}}
\put(13.00,54.00){\line(1,0){58.00}}
\put(10.00,54.00){\makebox(0,0)[cc]{$n_0$}}
\put(43.00,34.00){\makebox(0,0)[cc]{$s$}}
\put(32.00,45.00){\makebox(0,0)[cc]{$t'$}}
\put(50.00,70.00){\makebox(0,0)[cc]{$t$}}
\put(10.00,79.00){\makebox(0,0)[cc]{$N^*$}}
\put(80.00,80.00){\makebox(0,0)[cc]{$[t]\subseteq U_n$ in $r_2(\delta)$}}
\put(88.00,32.00){\makebox(0,0)[cc]
{$\swap(x_{\alpha_0},s,k)\cap U_n=\emptyset$ in $r_1(\delta)$}}
\put(95.00,48.00){\makebox(0,0)[cc]
{$\swap(x_{\alpha_0},t',N_1-1)\cap U_n=\emptyset$ in $r'(\delta) $}}
\end{picture}
\caption{The swap}
\end{figure}

\noindent This proves that
$((p',q_1\cup q_2),r_1(\delta)\cup r_2(\delta))\in \poset*\qposet$
for every $\delta$.

Finally we must show that
$$((p',q_1 \cup q_2),(r_1\cup r_2)\res\beta)\forces_\beta
\gamma\in Y_\beta$$
whenever
$(n,(\gamma,s,k))\in (r_1\cup r_2)(\beta)$ for some $n,s,k$.  But
by induction
$$((p',q_1\cup q_2),(r_1\cup r_2)\res\beta)$$ extends
both $((p_1,q_1),r_1\res\beta)$ and
$((p_2,q_2),r_2\res\beta)$, one of which does the required
forcing.

This proves the Claim.  The theorem now follows from the contradiction
that 
$$((p_1,q_1),r_1)\forces \tau\cap a_{\alpha_0}\subseteq N_1$$  
$$((p_2,q_2),r_2)\forces n_0\in\tau$$
where $n_0>N_1$ and
$$((p',q_1\cup q_2),r_1\cup r_2)\forces n_0\in a_{\alpha_0}$$
\qed

\bigskip

Remark.  The usual Q-set forcing 
kills the maximality of an almost disjoint family $X$. To see this
suppose $\{x_n:n<\omega\}\subseteq X$ and conditions 
are finite consistent sets of sentences of the form: ``$[s]\subseteq U_n$''
or ``$x\notin U_n$'' where $x\in X\setminus \{x_n:n<\omega\}$.  So
when we force we get a $G_\delta$ set so that
$$\cap_{n<\omega}U_n\cap X=\{x_n:n<\omega\}$$
In the generic extension we can find $\{k_n:n<\omega\}$ increasing
so that  
$$k_{n+1}\notin\cup_{i<n}x_i \rmand
\{y\subseteq \omega: k_{n+1}\in y\}\subseteq \bigcap_{i<n}U_i$$
Why? Given $p$ find $k_{n+1}>k_n$ not in any $x_i$ for $i<n$ or
in any $x$ mentioned in $p$ and put
$$p'=p\cup\{[s1]\subseteq U_i: i<n, s\in 2^{k_{n+1}-1}\}$$

\noindent 
But then $\{k_n:n<\omega\}$ is almost disjoint from all elements of $X$.

\bigskip
Remark.
Since there are perfect almost disjoint families, eg.,
  $$\aa=\{\{x\res n:n<\omega\}:x\in 2^\omega\}\subset P(2^{<\omega})$$
there are always MAD families of arbitrarily large Borel order.
Obviously a Q-set cannot have cardinality continuum, however 
a $\sigma$-set can.  

Define. $X\subseteq 2^\omega$ is a $\sigma$-set iff for every
Borel set $B\subseteq 2^\omega$ there exists a $G_\delta$ set
$G$ such that $B\cap X= G\cap X$.

A Sierpinski set is an example of $\sigma$-set
(Poprougenko, see Miller \cite{survey}).

\begin{theorem}
It is consistent with any cardinal arithmetic that there exists
a MAD $\sigma$-set of size the continuum.
\end{theorem}
\proof
This is an easy modification of the argument of the main theorem.  
Taking any countable transitive model $M$ first force a generic
MAD of size continuum, then do a finite support iteration of
length continuum to make it into a $\sigma$-set.
\qed

\bigskip 

Remark.  H.Woodin, see Larson \cite{larson}, has shown that 
if there exists a measurable Woodin cardinal
$\kappa$, and $V$ and $V[G]$ are both models of CH where $V[G]$ is a generic
extension using a partial order of size less than $\kappa$, then $V$ and $V[G]$
model exactly the same $\Sigma^2_1$ sentences.  The existence of a MAD
$\sigma$-set is a $\Sigma^2_1$ sentence.  It follows that

CH + there exists a measurable Woodin cardinal implies
there is a MAD $\sigma$-set.  

It is virtually certain that MAD $\sigma$-sets have nothing to do with
large cardinals, so we have the conjecture:

\begin{conjecture}
CH implies there exists a MAD $\sigma$-set.
\end{conjecture}

\begin{theorem}
The generic MAD set 
$\aa=\{a_\alpha:\alpha<\omega\}$ is
concentrated on $\{a_n:n<\omega\}$, ie. every open set containing
$\{a_n:n<\omega\}$ contains all but countably many elements of $A$.
\end{theorem}
\proof
Let $M$ be a countable standard model of ZFC and $G$ be 
$\poset$-generic over $M$.  Working in $M$ suppose 
$$\force \{a_n:n<\omega\}\subseteq U \mbox{ an open set}$$
Let $\Sigma\subseteq \poset$ be countable so that for
every $s\in 2^{<\omega}$ there exist a maximal antichain in
$\Sigma$ which decides ``$[s]\subseteq U$''.  

\bigskip\noindent
Claim. $\forces a_\alpha\in U$ for any $\alpha$ larger 
than any mentioned in $\Sigma$.

proof: Suppose not and let $(p,q)\force a_\alpha\notin U$.  Choose
some $n$ so that $n$ is not in the domain of $p$.  Let
$p'=p\cup \{(n,s)\}$ where $(\alpha,s)\in p$ and let 
$$q'=q\cup \bigcup \{(\{n,\beta\},k): (\{\alpha,\beta\},k)\in q\}$$
So $(p',q')\leq (p,q)$ and it says the same things about  $a_n$
and $a_\alpha$.  There exists $(\hat{p},\hat{q})\in\Sigma$ compatible
with $(p',q')$ such that 
$N_{\hat{p}}> N_{p}$ and 
$$(\hat{p},\hat{q})\forces [x_n\res N_{\hat{p}}]\subseteq U$$
Let $(p^*,q^*)$ extend both  $(p',q')$ and $(\hat{p},\hat{q})$.
Change $p^*$ to $r$ with same domain but $r(\alpha)=p^*(n)$
and other coordinate all the same.  But then
$(r,q'\cup \hat{q})$ is a common extension of both  
$(p',q')$ and $(\hat{p},\hat{q})$.  And this is a contradiction.

This proves the Claim and Theorem.
\qed
 
\begin{theorem}
CH implies exists a MAD family which is concentrated
on the finite subsets of $\omega$ and is a $\lambda$-set
(ie. every countable subset is a relative $G_\delta$. 
\end{theorem}
\proof
It is easy to construct a MAD family $\{a_\alpha:\alpha<\omega_1\}$ so
that if $f_\alpha:\omega\to a_\alpha$ is the strictly increasing
enumeration of $\alpha$, then for every
$\alpha<\beta$ we have that $f_\alpha<^* f_\beta$ and
for every $g\in\omega^\omega$ there exists $\alpha<\omega_1$
such that $g\leq^* f_\alpha$, ie. they form a scale.  Rothberger
(see Miller \cite{survey}) showed that any well-ordered subset
of $(\omega^\omega,\leq^*)$ is a $\lambda$-set and that any
$\omega_1$ -ordered unbounded set is concentrated on the
rationals.
\qed

\bigskip  
The same large cardinal results lead to the following
conjecture:

\begin{conjecture}
CH implies there exists a MAD family which is concentrated on
a countable subset of itself. 
\end{conjecture}

Paul Szeptycki pointed out that the Q-set forcing using in Theorem 1
can be used to prove the following:

\begin{theorem} It is relatively consistent that there exists
a Q-set $X\subseteq [\omega]^\omega$ satisfying the property
that for every $a\in [\omega]^\omega$ for all but countably
many $x\in X$ we have that $|x\cap a|=|x\setminus a|=\omega$,
ie $X$ is a strong splitting family.
\end{theorem}
\proof
We replace $\poset$ by the Cohen real partial order, i.e., just
drop the $q$'s from the $(p,q)$.  We use the same $\poset*\qposet$.
Note that in the basic argument for $p'$ we could have flipped
$p'(\alpha_0)(n_0)=1-p*(\alpha_0)(n_0)$ and $\alpha_0$ could be
any $\alpha<\omega_1$ not mentioned in $\Sigma$.
\qed

\begin{flushleft}
Arnold W. Miller \\
miller@math.wisc.edu \\
http://www.math.wisc.edu/$\sim$miller\\
University of Wisconsin-Madison \\
Department of Mathematics, Van Vleck Hall \\
480 Lincoln Drive \\
Madison, Wisconsin 53706-1388 \\
\end{flushleft}

\end{document}